%% file: High-Performance-Skew-Symmetric-EVP_PREPRINT_final.tex
\documentclass[5p,times]{elsarticle}

\usepackage{mymacros}

\usepackage{mathptmx}
\usepackage{amsmath}
\usepackage{amsfonts}
\usepackage{amssymb}

\usepackage{placeins}

\usepackage{algorithm}
\usepackage{algorithmicx}
\usepackage{algpseudocode}

\usepackage{subcaption}
 \usepackage{numprint}

\usepackage{pgfplots}
\pgfplotsset{compat=1.13}
\usetikzlibrary{shapes.multipart,decorations.markings,decorations.text, 
decorations.pathreplacing,patterns,positioning,intersections,calc,matrix}

\newtheorem{theorem}{Theorem}
\newtheorem{lemma}[theorem]{Lemma}

\usepackage{hyperref}

\journal{Parallel Computing}

\begin{document}

\begin{frontmatter}

\title{High Performance Solution of Skew-symmetric Eigenvalue Problems with Applications in Solving the Bethe-Salpeter Eigenvalue Problem}



\author[MPIMagdeburg]{Carolin Penke\corref{mycorrespondingauthor}}
\cortext[mycorrespondingauthor]{Corresponding author}
\ead{penke@mpi-magdeburg.mpg.de}

\author[MPCDF]{Andreas Marek}
\author[SOLgroup]{Christian Vorwerk}
\author[SOLgroup]{Claudia Draxl\fnref{BiGmax}}
\author[MPIMagdeburg]{Peter Benner\fnref{BiGmax}}
\fntext[BiGmax]{ These authors were supported by BiGmax, the Max Planck Society’s Research Network on Big-Data-Driven Materials Science.}

\address[MPIMagdeburg]{Computational Methods in Systems and
    Control Theory, Max Planck Institute for Dynamics of Complex Technical Systems, 
Germany}
\address[MPCDF]{Max Planck Computing and Data Facility, Garching,
Germany}
\address[SOLgroup]{Institut f\"ur Physik and IRIS Adlershof,
  Humboldt-Universit\"at zu Berlin, Berlin, Germany}

\begin{abstract}
We present a high-performance solver for dense skew-symmetric matrix eigenvalue
problems. Our work is motivated by applications in computational quantum
physics, where one solution approach to solve the Bethe-Salpeter equation involves
the solution of a large, dense, skew-symmetric eigenvalue problem. The
computed eigenpairs can be used to compute the optical absorption spectrum of
molecules and crystalline systems. One state-of-the art high-performance
solver package for symmetric matrices is the ELPA (Eigenvalue SoLvers for
Petascale Applications) library. We exploit a link between tridiagonal
skew-symmetric and symmetric matrices in order to extend the methods available
in ELPA to skew-symmetric matrices. This way, the presented solution method
can benefit from the optimizations available in ELPA that make it a
well-established, efficient and scalable library. The solution strategy is to
reduce a matrix to tridiagonal form, solve the tridiagonal eigenvalue problem
and perform a back-transformation for eigenvectors of interest. ELPA employs a
one-step or a two-step approach for the tridiagonalization of symmetric
matrices. We adapt these to suit the skew-symmetric case. The two-step
approach is generally faster as memory locality is exploited better. If all
eigenvectors are required, the performance improvement is counteracted by the
additional back transformation step. We exploit the symmetry in the spectrum
of skew-symmetric matrices, such that only half of the eigenpairs need to be
computed, making the two-step approach the favorable method. We compare
performance and scalability of our method to the only available
high-performance approach for skew-symmetric matrices, an indirect
route involving complex arithmetic. In total, we achieve a performance that is
up to 3.67 higher than the reference method using Intel's ScaLAPACK
implementation. Our method is freely available in the current release of the
ELPA library.
\end{abstract}


\end{frontmatter}


\section{Introduction}
A matrix $A\in\mathbb{R}^{n\times n}$ is called skew-symmetric when
$A=-A^{\tran}$, where $.^{\tran}$ denotes the transposition of a matrix. We are
interested in eigenvalues and eigenvectors of $A$. 

The symmetric eigenvalue problem, i.e. the case $A=A^{\tran}$, has been studied
in depth for many years. It lies at the core of many applications in different
areas such as electronic structure computations. Many methods for its solution
have been proposed \cite{GolV13} and successfully implemented. Optimized
libraries for many platforms are widely available
\cite{Lapack2,BlaCCetal97}. With the rise of more advanced computer
architectures and more powerful supercomputers, the solution of increasingly
complex problems comes within reach. Parallelizability and scalability become
key issues in algorithm development. The ELPA library \cite{MarBJetal14} is one
endeavor to tackle these challenges and provides highly competitive direct
solvers for symmetric (and Hermitian) eigenvalue problems running on distributed
memory machines such as compute clusters.

The skew-symmetric case \cite{WarG78} lacks the ubiquitous presence of its
symmetric counterpart and has not received the same extensive treatment. We
close this gap by extending the ELPA methodology to the skew-symmetric case. 

Our motivation stems from the connection to the Hamiltonian eigenvalue problem
which has many applications in control theory and model order reduction
\cite{BenKM05}. A real Hamiltonian matrix $H$ is connected to a symmetric matrix
$M$ via the matrix $J=\begin{bmatrix}0 & I\\-I&0\end{bmatrix}$, where $I$
denotes the identity matrix,
\begin{align*}
 M=JH.
\end{align*}
If $M$ is positive definite, in the following denoted by $M>0$, the Hamiltonian eigenvalue problem can be recast
into a skew-symmetric eigenvalue problem using the Cholesky factorization
$M=LL^{\tran}$. The eigenvalues of $H$ are given as eigenvalues of the
skew-symmetric matrix $L^{\tran}JL$ and eigenvectors can be transformed
accordingly.

This situation occurs for example in \cite{ShaJ16}, where a structure-preserving
method for the solution of the Bethe-Salpeter eigenvalue problem is described.
Solving the Bethe-Salpeter eigenvalue problem allows a prediction of optical
properties in condensed matter, a more accurate approach than currently used
ones, such as time-dependent density functional theory (TDDFT) \cite{Onida2002}.
In this application context, the condition $M > 0$ ultimately follows from much
weaker physical interactions represented in the off-diagonal
values~\cite{furcheDensityMatrixBased2001,
cizekStabilityConditionsSolutions1967}. When larger systems are of interest, the
resulting matrices easily become very high\--di\-men-sional. This calls for a
parallelizable and scalable algorithm. The solution of the corresponding
skew-symmetric eigenvalue problem can be accelerated via the developments
presented in this paper. 

The remaining paper is structured as follows. Section \ref{sec:sol} reintroduces
the methods used by ELPA and points out the necessary adaptations to make them
work for skew-symmetric matrices. The Bethe-Salpeter problem is presented in
Section \ref{Sec:BSE}. Section \ref{sec:numexp} provides performance results of
the ELPA extension, including GPU acceleration, and points out the speedup
achieved in the context of the Bethe-Salpeter eigenvalue problem.

\section{Solution Method} \label{sec:sol}
\subsection{Solving the Symmetric Eigenvalue Problem in ELPA}\label{sec:symsol}

The ELPA library \cite{MarBJetal14,AucBBetal11,AlvBBetal19} is a highly
optimized parallel MPI-based code \cite{For94}. It shows great scalability over
thousands of CPU cores and contains low-level optimizations targeting specific
compute architectures \cite{KusMKetal19}. When only a portion of eigenvalues and
eigenvectors are needed, this is exploited algorithmically and results in
performance benefits. We briefly describe the well-established procedure
employed by ELPA. This forms the basis of the method for skew-symmetric matrices
described in the next subsection.

ELPA contains functionality to deal with symmetric-definite generalized
eigenvalue problems. In this paper, we focus on the standard eigenvalue problem
for simplicity.  This is reasonable as it is the most common use case and forms
the basis of any method for generalized problems. We only consider real
skew-symmetric problems. The reason is that any skew-symmetric problem can be
transformed into a Hermitian eigenvalue problem by multiplying it with the
imaginary unit $i$. This problem can be solved using the available ELPA
functionality for complex matrices. For the real case this induces complex
arithmetic which should obviously be avoided, but for complex matrices this is a
viable approach.

We consider the symmetric eigenvalue problem, i.e. the orthogonal
diagonalization of a matrix,
\begin{align*}
Q^{\tran}AQ = \Lambda, 
\end{align*}
where $A=A^{\tran}\in\mathbb{R}^{n\times n}$ is the matrix whose eigenvalues are
sought. We are looking for the orthogonal eigenvector matrix $Q$ and the
diagonal matrix $\Lambda$ containing the eigenvalues. The solution is carried
out in the following steps.

\begin{enumerate}
 \item Reduce $A$ to tridiagonal form, i.e. find an orthogonal transformation
   $Q_{trd}$ s.t.
 \begin{align*}
  A_{trd} = Q_{trd}^{\tran}AQ_{trd}
 \end{align*}
 is tridiagonal. This is done by accumulating Householder transformations
 \begin{align*}
  Q_{trd} = Q_1Q_2\cdots Q_{n-1}, 
\end{align*}
  where $Q_i=I-\tau_i v_iv_i^{\tran}$ represents the $i$-th Householder
  transformation that reduces the $i$-th column and row of the updated
  $Q_{i-1}^{\tran}\cdots Q_1^{\tran}AQ_1\cdots Q_{i-1}$ to tridiagonal form.
  The matrices $Q_i$ are not formed explicitly but are represented by the
  Householder vectors $v_i$. These are stored in place of the eliminated
  columns of $A$.
\item Solve the tridiagonal eigenvalue problem, i.e. find orthogonal $Q_{diag}$
  s.t.
\begin{align*}
\Lambda= Q_{diag}^{\tran}A_{trd}Q_{diag}.
\end{align*}
In ELPA, this step employs a tridiagonal di\-vide-and-con\-quer scheme.
\item Transform the required eigenvectors back, i.e. perform the computation
\begin{align*}
 Q=Q_{trd}Q_{diag}.
\end{align*}
\end{enumerate}

The ELPA solver comes in two flavors which define the details of the
transformation steps, i.e Steps 1 and 3. ELPA1 works as described, the reduction
to tridiagonal form is performed in one step. ELPA2 splits the transformations
into two parts. Step 1 becomes
\begin{enumerate}
 \item 
 \begin{enumerate}
  \item Reduce $A$ to banded form, i.e. compute orthogonal $Q_{band}$ s.t.
  \begin{align*}
   A_{band}=Q_{band}^{\tran}AQ_{band}
  \end{align*}
  is a band matrix.
  \item Reduce the banded form to tridiagonal form, i.e. compute orthogonal
    $Q_{trd}$ s.t.
  \begin{align*}
   A_{trd}=Q_{trd}^{\tran}A_{band}Q_{trd}
  \end{align*}
  is tridiagonal.
 \end{enumerate}
\end{enumerate}
Accordingly, the back transformation step is split into two parts
\begin{enumerate}
\setcounter{enumi}{2}
 \item
 \begin{enumerate}
  \item Perform the back transformation corresponding to the band-to-tridiagonal
    reduction
  \begin{align*}
   \tilde{Q} = Q_{trd}Q_{diag}. 
  \end{align*}
  \item Perform the back transformation corresponding to the full-to-band reduction
  \begin{align*}
   Q = Q_{band}\tilde{Q}.
  \end{align*}
 \end{enumerate}
\end{enumerate}

The benefit of the two-step approach is that more efficient BLAS-3 procedures
can be used in the tridiagonalization process and an overlap of communication
and computation is possible. As a result, a lower runtime  can generally be
observed in the tridiagonalization, compared to the one-step approach. This
comes at the cost of more operations in the eigenvector back transformation due
to the extra step that has to be performed. Therefore, ELPA2 is superior to
ELPA1 in particular when only a portion of the eigenvectors is sought. In the
context of skew-symmetric eigenvalue problems, this becomes pivotal as the
purely imaginary eigenvalues come in pairs $\pm\lambda \iu$,
$\lambda\in\mathbb{R}$. The eigenvectors are given as the complex conjugates of
each other. It is therefore enough to compute half of the eigenvalues and
eigenvectors.

Both approaches are extended to skew-symmetric matrices in this work.

\subsection{Solving the Skew-symmetric Eigenvalue Problem}
Like a symmetric matrix, a skew-symmetric matrix can be reduced to tridiagonal
form using Householder transformations. A Householder transformation represents
a reflection onto a scaled first unit vector $e_1$. Let $H$ be a transformation
that acts on a vector $v$ s.t. $Hv=\alpha e_1$.  Obviously $-v$ is transformed
to $H(-v)=-\alpha e_1$ by the same $H$.  Therefore all tridiagonalization
methods that work on symmetric matrices, such as the ones implemented in ELPA,
can in principle work on skew-symmetric matrices as well. 

A skew-symmetric tridiagonal matrix is related to a symmetric one via the
following observation \cite{WarG78}.
\begin{lemma}\label{Lem:Symmetrization}
 With the unitary matrix $D=\text{diag}\{1,\iu,\iu^2,\dots,\iu^{n-1}\}$, where $\iu$ denotes the imaginary unit, $\alpha_j\in\mathbb{R}$, it holds
\begin{align}-\iu D^{\herm} {\scriptsize \begin{bmatrix}
		0 & \alpha_1 & & \\
		-\alpha_1 & 0 &\ddots & \\
		&\ddots &\ddots &\alpha_{n-1}\\
		& & -\alpha_{n-1} & 0                              
               \end{bmatrix}} D = {\scriptsize\begin{bmatrix}
		0 & \alpha_1 & & \\
		\alpha_1 & 0 &\ddots & \\
		&\ddots &\ddots &\alpha_{n-1}\\
		& & \alpha_{n-1} & 0                              
               \end{bmatrix}}.
 \end{align}
\end{lemma}
$.^{\herm}$ denotes the Hermitian transpose of a matrix. 

After the reduction to tridiagonal form, the symmetric tridiagonal system is solved using a divide-and-conquer method \cite{AucBBetal11}. As a first step of the back transformation, the resulting (real) eigenvectors have to be multiplied by the (complex) matrix $D$. Then the back transformations corresponding to the tridiagona\-lization take place. Algorithm \ref{alg:solution} outlines the process. It is very similar to the method employed for symmetric eigenvalue problems. The differences are the addition of step \ref{backtranD} and changes in the implementation, which are given in detail in Sections \ref{Impl:ELPA1} and \ref{Impl:ELPA2}.

\begin{algorithm}[t]
\caption{Solution of a Skew-symmetric Eigenvalue Problem} \label{alg:solution}
\label{dgeqrt}
\begin{algorithmic}[1]
\Require $A=-A^{\tran}\in\mathbb{R}^{n\times n}$
\Ensure Unitary eigenvectors $Q\in\mathbb{C}^{n\times n}$,
  $\lambda_1,\dots,\lambda_n \in\mathbb{R}$ s.t $Q^{\herm}AQ=
\text{diag}\{
\lambda_1\iu,\dots,\lambda_n\iu
\}.$

\State Reduce $A$ to tridiagonal form, i.e. generate $Q_{trd}$ s.t.\label{tran1}
\begin{gather*}
 Q_{trd}^{\tran}AQ_{trd} = A_{trd} ={\scriptsize\begin{bmatrix}
		0 & \alpha_1 & & \\
		-\alpha_1 & 0 &\ddots & \\
		&\ddots &\ddots &\alpha_{n-1}\\
		& & -\alpha_{n-1} & 0                              
               \end{bmatrix}} .
\end{gather*}

\State Solve the eigenvalue problem for the symmetric tridiagonal matrix $-\iu
  D^{\herm}A_{trd}D$, where $D=\text{diag}\{1,\iu,\iu^2,\dots,\iu^n\}$, i.e.
  generate $Q_{diag}$ s.t.
\begin{gather*}
 Q_{diag}^{\tran}{\scriptsize\begin{bmatrix}
		0 & \alpha_1 & & \\
		\alpha_1 & 0 &\ddots & \\
		&\ddots &\ddots &\alpha_{n-1}\\
		& & \alpha_{n-1} & 0                              
               \end{bmatrix}}Q_{diag}=
               {\scriptsize
               \begin{bmatrix}
               \lambda_1&&&\\
               &\lambda_2&&\\
               &&\ddots&\phantom{\alpha_{n-1}}\\
               &&\phantom{\alpha_{n-1}} &\lambda_n
               \end{bmatrix}}.
\end{gather*}

\State Back transformation corresponding to symmetrization (see Lemma \ref{Lem:Symmetrization}), i.e.
  compute\label{backtranD}
\begin{gather*}
 Q \leftarrow DQ_{diag}\in\mathbb{C}^{n\times n}. 
\end{gather*}
\State Back transformation corresponding to band-to-tridiagonal reduction, i.e.
compute \label{tran2}
\begin{gather*}
 Q\leftarrow Q_{trd}Q.
\end{gather*}

\end{algorithmic}
\end{algorithm}

In ELPA2 the transformation steps (\ref{tran1} and \ref{tran2} in Algorithm
\ref{alg:solution}) are both split into two parts as described in Section
\ref{sec:symsol}. 

\subsection{Implementation}

Extending ELPA for skew-symmetric matrices means adding the back transformation
step involving $D$. In contrast to symmetric matrices, skew-symmetric matrices
have complex eigenvectors and strictly imaginary eigenvalues. Computationally
complex values are introduced in Algorithm \ref{alg:solution} with $D$ in step
\ref{backtranD}. Further transformations have to be performed for the real and
the imaginary part individually. It is preferable to set up an array with
complex data type entries representing the eigenvectors as late as possible, so
that we can benefit from efficient routines in double precision.  The routines
for the eigenvector back transformation corresponding to tridiagonalization do
not change, because all they do is to apply Householder transformations to
non-symmetric (and non-skew-symmetric) matrices. They are applied on the real and imaginary part independently, realizing the complex back transformation in real arithmetic. The symmetric tridiagonal
eigensolver can be used as is. Ma\-king it aware of the zeros on the diagonal
might turn out to be numerically or computationally beneficial. 

We now examine the implementation of the two tridiagonalization approaches in
ELPA1 and ELPA2 in more detail. At many points in the original implementation,
symmetry of the matrix is assumed in order to avoid unnecessary computations and
to efficiently reuse data available in the cache. In this section we recollect
some details of the tridiagonal reduction in order to point out these instances.
Here, the implicit assumptions can be changed from ``symmetric'' to
``skew-symmetric'' by simple sign changes. 

ELPA is based on the well established and well documented 2D block-cyclic data
layout introduced by ScaLAPACK for load balancing reasons. It is therefore
compatible to ScaLAPACK and can act as a drop-in replacement while no ScaLAPACK
routines are used by ELPA itself.  In general, each process works on the part of
the matrix that was assigned to it. This chunk of data resides in the local
memory of the process. Communication between processes is realized via MPI. Each
process calls serial BLAS routines. Additional CUDA and OpenMP support is
available.

\subsubsection{Tridiagonalization in ELPA1}\label{Impl:ELPA1}
In ELPA1, the tridiagonalization is realized in one step using Householder
transformations. The computation of the Householder vectors is not affected by
the symmetry of a matrix. Essentially, the tridiagonalization of a matrix comes
down to a series of rank-2 updates \cite{MarRW68}, described in the following.
Given a Householder vector $v$, the update of the trailing submatrix is
performed as
\begin{align}\label{EQ:ELPA1:update1}
 A &\leftarrow (I-\tau v v^{\tran}) A (I-\tau v v^{\tran})\\
 &= A + v\underbrace{(0.5\tau^2 v^{\tran}Avv^{\tran} -\tau v^{\tran}A)}_{u_1^{\tran}} + \underbrace{(0.5\tau^2 vv^{\tran}Av- \tau A v )}_{u_2}v^{\tran}\label{EQ:ELPA1:update2}\\
 &= A + vu_1^{\tran} + u_2v^{\tran}\label{EQ:ELPA1:update3}\\
 &= A + \begin{bmatrix}
         v & u_2
        \end{bmatrix}
        \begin{bmatrix}
         u_1 & v
        \end{bmatrix}^{\tran}.\label{EQ:ELPA1:update4}
\end{align}

For symmetric matrices it holds $u_1=u_2$. This is assumed in the original ELPA
implementation. For skew-symmetric matrices it holds $u_1=-u_2$. In ELPA1, the
two matrices $\begin{bmatrix}
         v & u_2
        \end{bmatrix}$ and $\begin{bmatrix}
         u_1 & v
        \end{bmatrix}^{\tran}$ are stored explicitly. Actual updates are
        performed using \texttt{GEMM} and \texttt{GEMV} routines. The matrices
        differ in the data layout, i.e. which process owns which part of the
        matrix. After the vector $u_1$ is computed, it is transposed and
        redistributed to represent $u_2$ in $\begin{bmatrix}
         v & u_2
        \end{bmatrix}$. Here, for the skew-symmetric variant, a sign change is
        introduced. The skew-symmetric update now reads
\begin{align}
 A &\leftarrow A + \begin{bmatrix}
         v & -u_1
        \end{bmatrix}
        \begin{bmatrix}
         u_1 & v
        \end{bmatrix}^{\tran}.
\end{align}

During the computation of $u_1$, symmetry is assumed in the computation of
$A^{\tran}v$.  In particular, the code assumes that an off-diagonal matrix tile
is the same as in the transposed matrix. Another sign change corrects this
assumption for skew-symmetric matrices.

\subsubsection{Tridiagonalization in ELPA2}\label{Impl:ELPA2}
In ELPA2, the tridiagonalization is split into two parts. First, the matrix is
reduced to banded form, then to tridiagonal form. For the reduction to banded
form, the Householder vectors are computed by the process column owning the
diagonal block. They are accumulated in a triangular matrix
$T\in\mathbb{R}^{nb\times nb}$, where $nb$ is the block size.  The product of
Householder matrices is stored via its storage-efficient representation
\cite{SchV89}
\begin{gather}
 Q=H_1\cdots H_{nb} = I- VTV^{\tran},
\end{gather}
where $V=\begin{bmatrix}
          v_1&\cdots&v_{nb}
         \end{bmatrix}$ contains the Householder vectors. $H_i = I-\tau_i
         v_iv_i^{\tran}$ is the Householder matrix corresponding to the $i$-th
         Householder transformation.
         
In this context, the update of the matrix $A$ takes the following shape,
analogous to the direct tridiagonalization described in Section
\ref{Impl:ELPA1}.
\begin{align}
 A\leftarrow&\ (I-VTV^{\tran})^{\tran}A(I-VTV^{\tran})\\
 =&\ A + V\underbrace{(0.5T^{\tran}V^{\tran}AVTV^{\tran} -T^{\tran}V^{\tran}A) }_{U_1^{\tran}} \nonumber\\
 &+ \underbrace{(0.5 VT^{\tran}V^{\tran}AVT- AVT )}_{U_2}V^{\tran} \label{update}\\
 =&\ A + \begin{bmatrix}
         V & U_2
        \end{bmatrix}
        \begin{bmatrix}
         U_1 & V
        \end{bmatrix}^{\tran}. \label{updgemm}
\end{align}
It holds $U_1=U_2$ if $A$ is symmetric, and $U_1=-U_2$ if $A$ is skew-symmetric.
Each process computes the relevant parts of $U_1$ in a series of (serial) matrix
operations and updates the portion of $A$ that resides in its memory. Here, the
symmetry of $A$ is assumed and exploited at various points in the
implementation.  Sign changes have to be applied at these instances. 


For the banded-to-tridiagonal reduction, the matrix is redistributed in the form of
a 1D block cyclic data layout. Each process owns a diagonal and a subdiagonal
block. The reduction of a particular column introduces fill-in in the
neighboring block. The ``bulge-chasing'' is realized as a pipelined algorithm
where computation and communication can be overlapped by reordering certain
operations \cite{AucBBetal11,AucBH11}.

The update of the diagonal blocks takes the same form as in ELPA1 (Equations
\eqref{EQ:ELPA1:update1} to \eqref{EQ:ELPA1:update4}).
Here, no matrix multiplication is employed but BLAS-2 routines are used working
directly with the Householder vectors. It holds $u_1=u_2$ for symmetric $A$ and
$u_1=-u_2$ for skew-symmetric $A$. In the symmetric case, the update is realized
via a symmetric rank-2 update (\texttt{SYR2}). We implemented a skew-symmetric
variant of this routine which realizes the skew-symmetric rank-2 update $A
\leftarrow A - vu^{\tran} + uv^{\tran}$. For the setup of $u$, a skew-symmetric
variant of the BLAS routine performing a symmetric matrix vector product
(\texttt{SYMV}) is necessary.

The other parts of Algorithm \ref{alg:solution} are adopted from the symmetric
implementation without changes. The computation of Householder vectors, the
accumulation of the Householder transformations in a triangular matrix and the
update of the local block during reduction to banded form do not have to be
changed compared to symmetric ELPA. This is because they act on the lower part
of the matrix so that possible (skew-)symmetry has no effect. 

\section{The Bethe-Salpeter Eigenvalue Problem}\label{Sec:BSE}
\textit{Ab initio} spectroscopy aims to describe the excitations in condensed
matter from first principles, \textit{i.e.} without the input of any empirical
parameters. For light absorption and scattering, the Bethe-Salpeter Equation
(BSE) approach is the state-of-the-art methodology for both crystalline
systems\cite{Rohlfing1998,Benedict1998,Albrecht1998,SagA09,Onida2002} as well as
condensed molecular systems
\cite{Grossman2001,Faber2011,Cocchi2015,
Hirose2015}. This approach takes its name from the
\textit{Bethe-Salpeter Equation}~\cite{SalB51}, the equation of motion of the
electron-hole correlation function, as derived from many-body perturbation
theory~\cite{Strinati1988,Onida2002}. In practice, the problem of solving the
BSE is mapped to an effective eigenvalue problem. Specifically, its eigenvalues
and -states are employed to construct dielectric properties, such as the
spectral density, absorption spectrum, and the loss function~\cite{ShaJ16,Vorwerk2019}. An
appropriate discretization scheme leads to a finite-dimensional representation
in matrix form $H_{BS}$ that shows a particular block structure \cite{SaMK14}:
\begin{gather}\label{Eq:BSMatrix}
 H_{BS} = \begin{bmatrix}
           A& B\\
	   -\bar{B}& -\bar{A}
          \end{bmatrix} = 
          \begin{bmatrix}
           A & B \\
           -B^{\herm} & -A^{\tran}
          \end{bmatrix},\\          
          \quad A=A^{\herm},\quad B=B^{\tran} \in\mathbb{C}^{n\times n}.\nonumber
\end{gather}
Note that the Hermitian transpose $.^{\herm}$ as well as the regular transpose without complex conjugation $.^{\tran}$ play a role in this structure.

In general, we are interested in all eigenpairs of the Hamiltonian, as they
contain valuable information on the excitations of the system. Specifically,
they describe the \textit{bound excitons}, localized electron-hole
pairs that form due to correlation between an excited electron and a hole. The BSE
eigenstates are used to reconstruct the excitonic wavefunction and obtain the
excitonic binding energy.

In this paper, we present a solution strategy for the most general formulation
of the BSE problem. As such, $A$ and $B$ are generally dense and complex-valued,
which holds in the case of excitations in condensed matter.

$H_{BS}$ belongs to the slightly more general class of $J$-sym\-met\-ric
matrices \cite{BenFY18}. This class of matrices display a symmetry
$(\lambda,-\lambda)$ in the spectrum. The additional structure in $H_{BS}$ leads
to an additional symmetry $(\lambda,  -\lambda, \bar{\lambda}, -\bar{\lambda})$
and a relation between the corresponding eigenvectors. Following \cite{ShaJ16},
we consider the definite Bethe-Salpeter eigenvalue problem. $H_{BS}$ is called
definite when the property
\begin{gather}\label{BSdefinite}\begin{bmatrix} I_n & 0\\
	0 & -I_n\end{bmatrix} H_{BS} = 
	\begin{bmatrix}
	A & B\\
	\bar{B} & \bar{A}
	\end{bmatrix}	> 0
\end{gather}
is fulfilled, which often holds in practice. In this case, the eigenvalues are
real and therefore come in pairs $(\lambda, - \lambda)$. The method presented in this work relies on this assumption. 

We aim for a solution method that preserves this structure under the influence
of inevitable numerical errors, i.e. that guarantees that the eigenvalues come
in pairs or quadruples, respectively. General methods for eigenvalue problems,
such as the QR/QZ algorithm, destroy this property. In this case it is not clear
anymore which eigenpairs correspond to the same excitation state. 

A structure-preserving method running in parallel on distributed memory systems
is developed in \cite{ShaJ16} and has been made available as BSEPACK. It relies
on assumption \eqref{BSdefinite} and exploits a connection to a Hamiltonian
eigenvalue problem  given in the following Theorem.

\begin{theorem}\label{Thm:HamiltonianTransform}
Let 
  $Q=\frac{1}{\sqrt{2}}\begin{bmatrix}
                       I & -\iu I\\
                       I & \iu I
                      \end{bmatrix}$,
 then $Q$ is unitary and
 \begin{align*}
  Q^{\herm} \begin{bmatrix}
       A & B \\ 
       -\bar{B}& -\bar{A}
      \end{bmatrix}Q  = i
      \begin{bmatrix}
       \Imag{A+B} & -\Real{A-B} \\
       \Real{A+B} & \Imag{A-B}
      \end{bmatrix}  =: \iu H,
 \end{align*}
where $H$ is real Hamiltonian, i.e.  $JH = (JH)^{\tran}$ with\\ $J=\begin{bmatrix}
                                                             0 & I\\
                                                             -I & 0
                                                            \end{bmatrix}$.
\end{theorem}

Let 
\begin{align}\label{Eq:SymmetricM}
M=JH =\begin{bmatrix}
       \Real{A+B} & \Imag{A-B}\\
       -\Imag{A+B} & \Real{A-B}
      \end{bmatrix} 
\end{align}
be the symmetric matrix associated with the Hamiltonian matrix $H$. Its positive
definiteness follows from property \eqref{BSdefinite}, which can be seen in the following way. 
Let the matrices $S$ and $\Omega$ be given as
\begin{align}
 S= \begin{bmatrix}
     I_n&\\
     &-I_n
    \end{bmatrix},\qquad \Omega = 
    \begin{bmatrix}
     A&B\\
     \bar{B} & \bar{A}
    \end{bmatrix},
\end{align}
i.e. $H_{BS}=S\Omega$. With the matrix $Q$ from Theorem \ref{Thm:HamiltonianTransform} we have
\begin{align}\label{Eq:Mbuild}
 M = -\iu J Q^{\herm}S\Omega Q.
\end{align}
It is easily verified that
\begin{align}
 -\iu J Q^{\herm}SQ=I_n,
\end{align}
i.e. $-\iu J Q^{\herm}S$ is the inverse of $Q$. The construction of $M$ \eqref{Eq:Mbuild} can therefore be seen as a similarity transformation of $\Omega$. If $\Omega$ is positive definite \eqref{BSdefinite}, so is M. The method described in
\cite{ShaJ16} relies on this property in order to guarantee the existence of the Cholesky factorization of $M$.

It performs the following steps.
\begin{enumerate}
 \item Construct $M$ as in \eqref{Eq:SymmetricM}.
 \item Compute a Cholesky factorization $M=LL^{\tran}$.
 \item Compute eigenpairs of the skew-symmetric matrix $L^{\tran}JL$, where $J=\begin{bmatrix}
                                                             0 & I\\
                                                             -I & 0
                                                            \end{bmatrix}$.
 \item Perform the eigenvector back transformation associated with Cholesky
   factorization and transformation to Hamiltonian form (Theorem
    \ref{Thm:HamiltonianTransform}).
\end{enumerate}

The eigenvalues and eigenvectors can be used to compute the optical absorption
spectrum of the material in a postprocessing step. 

The main workload is given as the solution of a skew-sym\-met\-ric eigenvalue
problem (Step 3). As a proof of concept, solution routines for the symmetric
eigenvalue problem from the ScaLAPACK reference implementation
\cite{BlaCCetal97} were adapted to the skew-symmetric setting. The matrix is
reduced to tridiagonal form using Householder transformations. The tridiagonal
eigenvalue problem is solved via bisection and inverse iteration. 

The ScaLAPACK reference implementation is not regarded as a state-of-the art solver
library. When performance and scalability are issues, one generally turns to
professionally maintained and optimized libraries such as ELPA
\cite{MarBJetal14} or vendor-specific implementations such as Intel's MKL.
Within BSE\-PACK, ScaLAPACK can be substituted by ELPA working on skew-symmetric
matrices. The resulting performance benefits are discussed in Section
\ref{Sec:NumBSE}.

\section{Numerical Experiments}\label{sec:numexp}
\subsection{ELPA Benchmarks}
\newlength\figureheight
\newlength\figurewidth 
\setlength\figureheight{0.26\textwidth}
\setlength\figurewidth{0.4\textwidth}
\begin{figure}[t]
\input{ELPASS_F1.tex}
 \caption{Scaling of the ELPA solver for skew-symmetric matrices. For comparison
  the runtimes for the alternative solution method via complex Hermitian solvers
  is included. Here, ELPA and Intel's MKL 2018 routines \texttt{pzheevd} and
  \texttt{pzheevr} are used. The matrix has a size of $n=20\,000$.
  \label{fig:scalability}}
\end{figure}
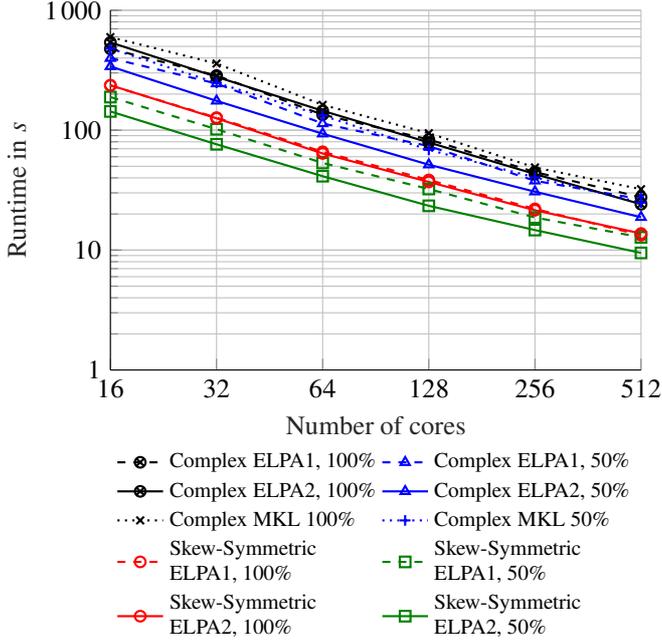

In this section we present performance results for the skew-symmetric ELPA
extension. All test programs are run on the \textit{mechthild} compute cluster,
located at the Max Planck Institute for Dynamics of Complex Technical Systems in
Magdeburg, Germany. Up to 32 nodes are used, which consist of 2 Intel Xeon
Silver 4110 (Skylake) processors with 8 cores each, running at 2.1 GHz. The
Intel compiler, MPI library and MKL in the 2018 version are used in all test
programs.  The computations use randomly generated skew-symmetric matrices in
double-precision. 

Figure \ref{fig:scalability} shows the resulting performance and the scaling
properties of ELPA for a medium sized skew-symmetric matrix ($n=20\,000$). As an
alternative to the approach described in this work, the skew-symmetric matrix
can be multiplied with  the imaginary unit $i$. The resulting complex Hermitian
matrix can be diagonalized using available methods in ELPA or Intel's ScaLAPACK
implementation shipped with the MKL. This represents the only previously available approach
to solve skew-symmetric eigenvalue problems in a massively parallel
high-performance setting. 

For skew-symmetric matrices, only 50\% of eigenvalues and eigenvectors need to
be computed, as they are purely imaginary and come in pairs $\pm\lambda \iu,
\lambda\in\mathbb{R}$. The runtime measurements for 100\%  are included for
reference. 
\begin{table}[t]
\center
\caption{Execution time speedups achieved by different aspects of the solution approach.}
  \label{tab:Scalability}
\begin{tabular}{|p{0.0375\textwidth}||>{\centering\arraybackslash}p{0.07\textwidth}||>{\centering\arraybackslash}p{0.07\textwidth}|>{\centering\arraybackslash}p{0.07\textwidth}||>{\centering\arraybackslash}p{0.07\textwidth}|}\hline
\scriptsize\#Cores &\scriptsize Compl. ELPA2 100\% vs. Compl. MKL 100 \% &\scriptsize Compl. ELPA2 50\% vs. Compl. MKL 50\% &\scriptsize Skew-Sym. ELPA2 50\% vs. Compl. ELPA2 50\% &\scriptsize Skew-Sym. ELPA2 50\% vs. Compl. MKL 50\%      \\\hline
16& 1.10&1.41&2.33&3.28\\
32& 1.29&1.41&2.30&3.24\\
64& 1.11&1.40&2.32&3.25\\
128&1.18&1.33&2.20&2.93\\
256&1.17&1.28&2.16&2.76\\
512&1.21&1.51&1.87&2.82\\\hline
\end{tabular}
\end{table}

Figure \ref{fig:scalability} shows that all approaches display good scalability
in the examined setting. Skew-symmetric ELPA runs 2.76 to 3.28 times faster than
the complex MKL based solver, where both only compute 50\% of eigenpairs. The
data gives further insight into how this improvement is achieved. Table
\ref{tab:Scalability} compares the runtimes for different solvers and presents
the achieved speed\-ups.  When we compare complex 100\% solvers, ELPA already
improves performance by a factor of $1.1$ to $1.29$ (column 2 in Table
\ref{tab:Scalability}).  When all eigenpairs are computed, EL\-PA1 and ELPA2
yield very similar runtime results which is why only ELPA2 is considered in
Table \ref{tab:Scalability}. The two-step approach employed by ELPA2 pays off in
particular when not all eigenpairs are sought, which is the case here. When
complex 50\% solvers are compared (ELPA2 vs. MKL, column 3 in Table
\ref{tab:Scalability}), the achieved speedup increases to a value between $1.28$
and $1.51$. The largest impact on the performance is caused by avoiding complex
arithmetic. This is represented by the speedup achieved by the skew-symmetric 50\% ELPA2 implementation compared to the complex 50\% ELPA2 implementation
(column 4 of Table \ref{tab:Scalability}). This accounts for an additional
speedup of 1.87 to 2.33. 

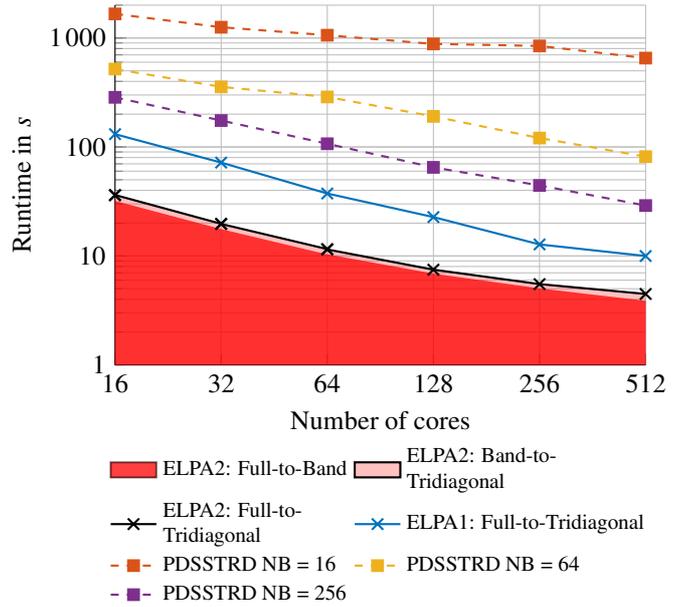
\begin{figure}[t]
\centering
\input{ELPASS_F2.tex}
 \caption{Scaling of the tridiagonalization in two steps (ELPA2) and one step
  (ELPA1). We compare it to the runtimes of the tridiagonalization routine for
  skew-symmetric matrices \texttt{PDSSTRD} available in BSEPACK \cite{ShaJ16}
  for different block sizes $NB$. The matrix size is
  $n=20\,000$.\label{fig:scalabilityTRD}}
\end{figure}

The tridiagonalization is an essential step in every con\-si\-der\-ed solution
scheme and contributes a significant portion of the execution time. The fewer
eigenpairs are sought, the more dominant it becomes with respect to computation
time. Figure \ref{fig:scalabilityTRD} displays the runtimes and scalability of
available tridiagonalization techniques for skew-symmetric matrices. As an
alternative implementation to the presented approaches there is a
tridiagonalization routine \texttt{PDSSTRD} shipped in BSEPACK \cite{ShaJ16}. It
is an adapted version of the ScaLAPACK reference implementation.

All discussed implementations are based on the 2D-block-cyclic data distribution
established by ScaLAPACK. Here, the matrix is divided into blocks of a certain
size $NB$. The blocks are distributed to processes organized in a 2D grid in a
cyclic manner. Typically, the block size is a parameter chosen once in a
software project. The data redistribution to data layouts defined by other block
sizes is avoided as this involves expensive all-to-all communication.  The main
disadvantage of the \texttt{PDSSTRD} routine is that it is very susceptible to
the chosen block size, both with regard to scalability and overall performance.
This makes it less suitable to be included in larger software projects, where
the block size is a parameter predefined by other factors. ELPA (both the one
and two-step version) on the other hand does not have this problem and performs
equally well for all data layouts \cite{BenMP18}.

Figure \ref{fig:scalabilityTRD} also displays the advantage of the two-step
tridiagonalization over the one-step approach. Here the performance is dominated
by the first step, i.e. the reduction to banded form.
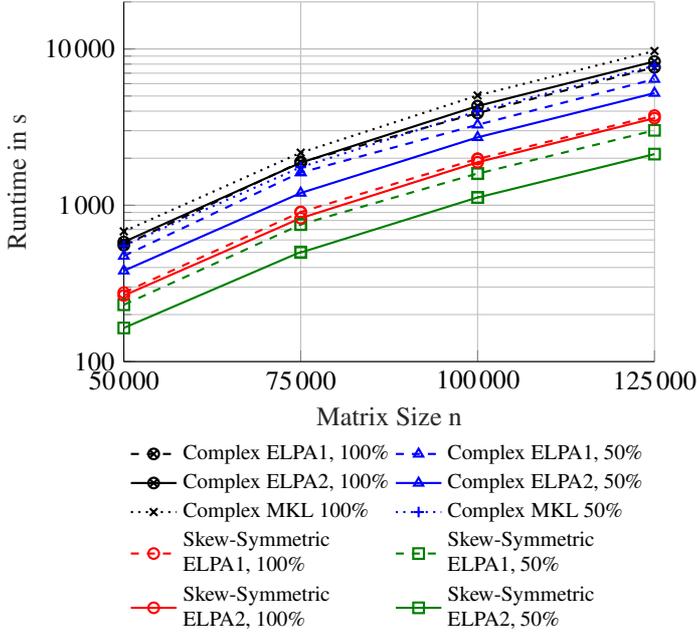
\begin{figure}[t]
\centering
\input{ELPASS_F3.tex}
 \caption{Runtimes for solving eigenvalue problems of larger
  sizes.\label{fig:largematrices} 256 CPU cores were used, i.e. 16 nodes on the
  \textit{mechthild} compute cluster.}
\end{figure}
In the context of electronic structure computations, the matrices of interest
can become extremely large. Figure \ref{fig:largematrices} displays the achieved
runtime improvements for larger matrices up to a size of $n=125\,000$. The
individual speedups are presented in Table \ref{tab:largematrices}. For large
matrices we achieve a speedup of up to 3.67 compared to the available MKL
routine. 

\begin{table}[t]
\center
\caption{Execution time speedups achieved by different aspects of the solution
  approach.}
  \label{tab:largematrices}
\begin{tabular}{|p{0.06\textwidth}||>{\centering\arraybackslash}p{0.07\textwidth}||>{\centering\arraybackslash}p{0.07\textwidth}|>{\centering\arraybackslash}p{0.07\textwidth}||>{\centering\arraybackslash}p{0.07\textwidth}|}\hline
Matrix size &\scriptsize Compl. ELPA2 100\% vs. Compl. MKL 100 \% &\scriptsize
  Compl. ELPA2 50\% vs. Compl. MKL 50\% &\scriptsize Skew-Sym. ELPA2 50\% vs.
  Compl. ELPA2 50\% &\scriptsize Skew-Sym. ELPA2 50\% vs. Compl. MKL 50\%
  \\\hline
50\,000&   1.17&1.45&2.32& 3.35\\
75\,000&   1.16&1.46&2.39& 3.50\\
100\,000&   1.17&1.47&2.42&3.57\\
125\,000&1.17&1.49&2.46&   3.67\\\hline
\end{tabular}
\end{table}

\subsubsection{GPU Acceleration}
\begin{figure}[t]
\centering
\input{1NodeGPU.tex}
 \caption{\label{Fig:GPU} Runtimes for solving eigenvalue problems on one node
 on the \textit{mechthild} compute cluster employing a GPU.}
\end{figure}
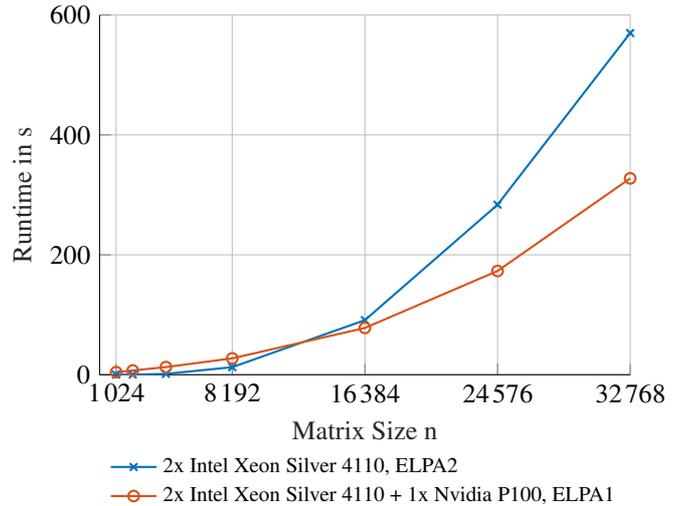
For the 1-step tridiagonalization approach (ELPA1), there is a GPU-accelerated
version available that gets shipped with the ELPA library \cite{KusLM19}. The
design approach is to stick with the same code base as the CPU-only version, and
offload compute-intense parts, such as BLAS-3 operations, to the GPU in order to
benefit from its massive parallelism. This is done using the CUBLAS library
provided by NVIDIA.  Because ELPA2 employs more fine-grained communication
patterns, this approach works best for ELPA1. Here, the performance can benefit
when the computational intensity is high enough, i.e. when big chunks of data
are being worked on by the GPU. 

Figure \ref{Fig:GPU} shows the performance that can be achieved on one node of
the \texttt{mechthild} compute cluster, that is equipped with an NVIDIA P100 GPU
as an accelerator device. The GPU version is based on ELPA1 and therefore does
not benefit from the faster tridiagonalization in ELPA2 (see Figure
\ref{fig:scalabilityTRD} and the discussion in the previous section). Despite
this fact, the GPU-accelerated ELPA1 version eventually outperforms the ELPA2
CPU-only version, if the matrix is large enough. In our case the turning point
is at around $n=15\,000$. For smaller matrices the additional work of setting up
the CUDA environment and transferring the matrix counteracts any possible
performance benefits and results in a larger runtime. For matrices of size
$n=32\,768$ employing the GPU can reduce the runtime from 570 seconds to 328
seconds, i.e. by 41\%. 

The take-away message of these results is the following. If nodes equipped with
GPUs are available and to be utilized, it is important to make sure each node
has enough data to work on. This way, the available resources are used most
efficiently.

\subsection{Accelerating BSEPACK}\label{Sec:NumBSE}
\setlength\figureheight{0.26\textwidth}
\setlength\figurewidth{0.4\textwidth}
\begin{figure}[t]
\input{AccBSE.tex}
 \caption{Scaling of the direct, complex BSE\-PACK eigenvalue solver for
  computing the optical absorption spectrum of hexagonal boron nitride. The
  Bethe-Salpeter matrix \eqref{Eq:BSMatrix} has a size of
  $51\,200$.\label{Fig:BSEscalability}}
\end{figure}
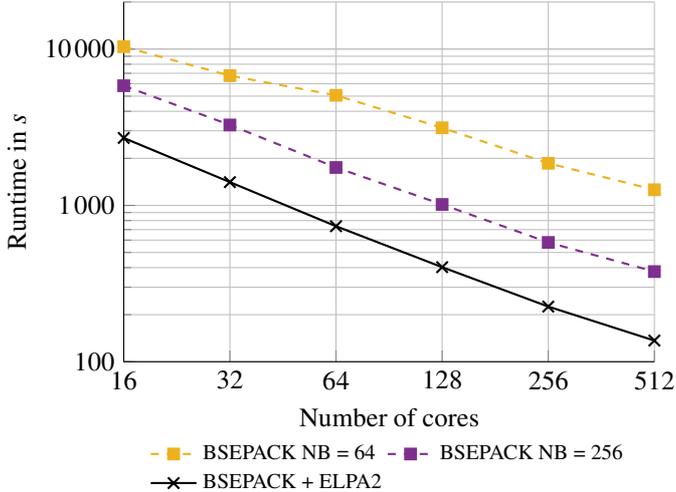

We consider the performance improvements that can be a\-chieved by using the
newly developed skew-symmetric eigenvalue solver in the BSEPACK \cite{ShaJ16}
software, described in Section \ref{Sec:BSE}. In this procedure, Step
3, the computation of eigenpairs of the skew-symmetric matrix $L^{\tran}JL$, is
now performed by the ELPA library. 

To demonstrate the speedup, we consider the example of hexagonal
boron nitride at a fixed size of the BSE Hamiltonian. The excitations in
hexagonal boron nitride are widely studied both experimentally and theoretically
\cite{Cappellini2001,Blase1995,Galambose2011,Fugallo2015,Cudazzo2016,
Koskelo2017,Aggoune2018},
as its wide band gap and the layered geometrical structure yield strong effects
of electron-hole correlation, such as the formation of bound excitons.
Previous studies have shown that the BSE approach yields the optical absorption
and excitonic properties with high accuracy. In our calculations, the BSE
Hamiltonian is constructed on a $16\times 16 \times 4$ $\mathbf{k}$-grid in the
1st Brillouin zone, the 5 highest valence and 5 lowest conduction bands are
employed to construct the transition space, leading to a matrix size of $2
\times 16 \times 16 \times 4 \times 5 \times 5=51200$. In the calculation of
the BSE Hamiltonian, single-particle wavefunctions and the static dielectric
function are expanded in plane waves with a cut-off of 387 eV and 132 eV,
respectively. The static dielectric function is obtained from ABINIT
\cite{ABINIT2016}, while the
BSE Hamiltonian is constructed using the EXC code \cite{EXC2020}.

Figure \ref{Fig:BSEscalability} displays the achieved runtimes of BSEPACK for
this fixed-size matrix for different core counts. We compare the original
version and a version that employs ELPA. The performance of the original solver
is highly dependent on the chosen block size (see also Figure
\ref{fig:scalabilityTRD}). This parameter determines how the matrix is
distributed to the available processes in the form of a 2D block-cyclic data layout.
The default is given as $NB=64$, but choosing a larger block size can increase
the performance dramatically, as can be seen in Figure \ref{Fig:BSEscalability}
for $NB=256$. Typically, software packages (e.g.
\cite{Gulans2014,Vorwerk2019}) developed
for electronic structure computations are large and contain many features,
implementing methods for different quantities of interest. The block size is
typically predetermined by other considerations. It would mean a serious effort
to change it, in order to optimize just one building block of the software.
Furthermore the optimal block size of the original BSEPACK is probably dependent
on the given hardware and the given matrix size. Autotuning frameworks could
help, but are also very costly and impose an additional implementation effort. A
software, that does not show this kind of runtime dependency is greatly
preferable. Employing ELPA for the main computational task in BSE\-PACK fulfills
this requirement. The performance of ELPA is independent of the chosen $NB$,
because the block size on the node level for optimal cache use is decoupled from
the block size defining the multi-node data layout. 

The ELPA-accelerated version is up to 9.22 times as fast as the original code
with the default block size.  Even when the block size is increased, using the
new solver always yields a better performance. In the case of $NB=256$, the
ELPA-version still performs up to 2.76 times as fast. Choosing even larger block
sizes has in general no further positive effect on the performance of the
original BSEPACK. Employing ELPA also leads to an improved scalability over the
number of cores.

\section{Conclusions}

We have presented a strategy to extend existing solver libraries for symmetric
eigenvalue problems to the skew-sym\-met\-ric case. Applying these ideas to the
ELPA library, makes it possible to compute eigenvalues and eigenvectors of large
skew-symmetric matrices in parallel with a high level of efficiency and
scalability. We benefit from the maturity of the ELPA software project, where
many optimizations have been realized over the years. All of these, including
GPU support, find their way into the presented skew-symmetric solver. As far as
we know, no other solvers dedicated to the skew-symmetric eigenvalue problem
exist in an HPC setting. It is always possible to solve a complex Hermitian
eigenvalue problem instead of a skew-symmetric one. Our newly developed solver
outperforms this strategy, implemented via Intel MKL ScaLAPACK, by a factor of
3. We also observe an increase in performance concerning the Bethe-Salpeter
eigenvalue problem. Here we improve the runtime of available routines by a
factor of almost 10, making the BSEPACK library with ELPA a viable choice as a
building block for larger electronic structure packages.

\section{Acknowledgment}
We thank Francesco Sottile for fruitful discussion and his support in generating
the BSE Hamiltonian for hexagonal BN.


\input{High-Performance-Skew-Symmetric-EVP.bbl}
\end{document}

%% file: ELPASS_F1.tex
%
%
\begin{tikzpicture}

\begin{axis}[%
/pgf/number format/.cd,
        1000 sep={\,},
log basis x=16,
    log ticks with fixed point,
width=0.951\figurewidth,
height=\figureheight,
at={(0\figurewidth,0\figureheight)},
scale only axis,
xmode=log,
xmin=16,
xmax=512,
xtick={ 16,  32,  64, 128, 256, 512},
xminorticks=true,
xlabel style={font=\color{white!15!black}},
xlabel={Number of cores},
ymode=log,
ymin=1,
ymax=1000,
yminorticks=true,
ylabel style={font=\color{white!15!black}},
ylabel={Runtime in $s$},
axis background/.style={fill=white},
axis x line*=bottom,
axis y line*=left,
xmajorgrids,
xminorgrids,
ymajorgrids,
yminorgrids,
legend columns=2,
legend style={at={(0.5,-.2)}, anchor=north, legend cell align=left, align=left, draw=white!15!black, font=\footnotesize, draw=none},
every axis plot/.append style={thick}
] 
\addplot [color=black, dashed, mark=otimes, mark options={solid, black}]
  table[row sep=crcr]{%
16	476.739407\\
32	285.529594\\
64	134.605456\\
128	83.666212\\
256	44.737435\\
512	27.566232\\
};
\addlegendentry{Complex ELPA1, 100\%}

\addplot [color=blue, dashed, mark=triangle, mark options={solid, blue}]
  table[row sep=crcr]{%
16	398.50685\\
32	244.713447\\
64	113.620397\\
128	73.665949\\
256	37.829783\\
512	26.529185\\
};
\addlegendentry{Complex ELPA1, 50\%}

\addplot [color=black, mark=otimes, mark options={solid, black}]
  table[row sep=crcr]{%
16	539.732537\\
32	278.512065\\
64	145.034454\\
128	78.977783\\
256	43.07808\\
512	24.113866\\
};
\addlegendentry{Complex ELPA2, 100\%}

\addplot [color=blue, mark=triangle, mark options={solid, blue}]
  table[row sep=crcr]{%
16	340.047587\\
32	175.580578\\
64	93.4223\\
128	51.525228\\
256	30.686918\\
512	18.79392\\
};
\addlegendentry{Complex ELPA2, 50\%}

\addplot [color=black, dotted, mark=x, mark options={solid, black}]
  table[row sep=crcr]{%
16	597.416659\\
32	359.667624\\
64	162.606278\\
128	94.324431\\
256	48.902167\\
512	32.195974\\
};
\addlegendentry{Complex MKL 100\%}

\addplot [color=blue, dotted, mark=+, mark options={solid, blue}]
  table[row sep=crcr]{%
16	478.772203\\
32	248.303823\\
64	132.294128\\
128	68.556339\\
256	40.357521\\
512	24.893169\\
};
\addlegendentry{Complex MKL 50\%}

\addplot [color=red, dashed, mark=o, mark options={solid, red}]
  table[row sep=crcr]{%
16	235.756207\\
32	127.212542\\
64	65.980196\\
128	38.462357\\
256	22.052173\\
512	13.290403\\
};
\addlegendentry{Skew-Symmetric\\ ELPA1, 100\%}

\addplot [color=black!50!green, dashed, mark=square, mark options={solid, black!50!green}]
  table[row sep=crcr]{%
16	188.837467\\
32	101.841065\\
64	53.274354\\
128	32.303058\\
256	18.704349\\
512	12.767152\\
};
\addlegendentry{Skew-Symmetric \\ ELPA1, 50\%}

\addplot [color=red, mark=o, mark options={solid, red}]
  table[row sep=crcr]{%
16	236.992734\\
32	125.521974\\
64	63.928667\\
128	37.007031\\
256	21.489789\\
512	13.747303\\
};
\addlegendentry{Skew-Symmetric\\ ELPA2, 100\%}
%
%
%
%
%

\addplot [color=black!50!green, mark=square, mark options={solid, black!50!green}]
  table[row sep=crcr]{%
16	143.28074\\
32	76.244042\\
64	41.31252\\
128	23.39518\\
256	14.687848\\
512	9.467844\\
};
\addlegendentry{Skew-Symmetric\\ ELPA2, 50\%}

\end{axis}
\end{tikzpicture}%

%% file: ELPASS_F2.tex
%
%
\definecolor{mycolor1}{rgb}{0.00000,0.44700,0.74100}%
\definecolor{mycolor2}{rgb}{0.85000,0.32500,0.09800}%
\definecolor{mycolor3}{rgb}{0.92900,0.69400,0.12500}%
\definecolor{mycolor4}{rgb}{0.49400,0.18400,0.55600}%
\begin{tikzpicture}

\begin{axis}[%
/pgf/number format/.cd,
        1000 sep={\,},
log basis x=16,
log ticks with fixed point,
xtick={ 16,  32,  64, 128, 256, 512},
xlabel={Number of cores},
width=0.951\figurewidth,
height=\figureheight,
at={(0\figurewidth,0\figureheight)},
scale only axis,
xmode=log,
xmin=16,
xmax=512,
ylabel={Runtime in $s$},
xminorticks=true,
ymode=log,
ymin=1,
ymax=2000,
yminorticks=true,
axis background/.style={fill=white},
axis x line*=bottom,
axis y line*=left,
xmajorgrids,
xminorgrids,
ymajorgrids,
yminorgrids,
legend columns=2,
legend style={at={(0.5,-.2)}, anchor=north, legend cell align=left, align=left, draw=none, font=\footnotesize},
every axis plot/.append style={thick}
]
\addplot[fill=red, fill opacity=0.75, draw=none, draw opacity=0.5, area legend] table[row sep=crcr, stack plots=y]{%
16	31.751869\\
16	32.015465\\
32	17.64709\\
32	17.644777\\
64	10.316835\\
64	10.357729\\
128	6.886264\\
256	5.0385115\\
512	3.879271\\
}
\closedcycle;
\addlegendentry{ELPA2: Full-to-Band}

\addplot[fill=red, fill opacity=0.25, draw=none, area legend] table[row sep=crcr, stack plots=y]{%
16	4.2402\\
16	4.265257\\
32	2.068185\\
32	2.019592\\
64	1.242669\\
64	1.109548\\
128	0.596949\\
256	0.4761505\\
512	0.6003145\\
}
\closedcycle;
\addlegendentry{ELPA2: Band-to-\\Tridiagonal}

\addplot[draw=black, mark=x, mark size=3pt] table[row sep=crcr, stack plots=y]{%
16	0.000000000001 \\
16	0.000000000001 \\
32	0.000000000001 \\
32	0.000000000001 \\
64	0.000000000001 \\
64	0.000000000001 \\
128	0.000000000001 \\
256	0.000000000001 \\
512	0.000000000001 \\
};
\addlegendentry{ELPA2: Full-to-\\Tridiagonal}

\addplot [color=mycolor1,  mark=x, mark options={solid, mycolor1}, mark size=3pt, stack plots=false]
  table[row sep=crcr]{%
16	131.085559\\
32	71.8952235\\
64	37.4036185\\
128	22.7285615\\
256	12.7525015\\
512	9.9853605\\
};
\addlegendentry{ELPA1: Full-to-Tridiagonal}

\addplot [color=mycolor2, dashed, mark=square*, mark options={solid, mycolor2, fill}, stack plots=false]
  table[row sep=crcr]{%
16	1663.857546\\
32	1255.269306\\
64	1060.342041\\
128	882.194117\\
256	845.680006\\
512	654.20812\\
};
\addlegendentry{PDSSTRD NB = 16}
\addplot [color=mycolor3, dashed, mark=square*, mark options={solid, mycolor3, fill}]
  table[row sep=crcr, stack plots=false]{%
16	518.885653\\
32	356.056553\\
64	287.77952\\
128	190.88687\\
256	120.752605\\
512	81.502387\\
};
\addlegendentry{PDSSTRD NB = 64}

\addplot [color=mycolor4, dashed, mark=square*, mark options={solid, mycolor4, fill}]
  table[row sep=crcr, stack plots=false]{%
16	285.450451\\
32	174.786484\\
64	107.028445\\
128	64.996395\\
256	44.357952\\
512	29.006982\\
};
\addlegendentry{PDSSTRD NB = 256}

\end{axis}
\end{tikzpicture}%

%% file: ELPASS_F3.tex
%
%
\begin{tikzpicture}

\begin{axis}[%
/pgf/number format/.cd,
        1000 sep={\,},
scaled x ticks = false,
x tick label style = {/pgf/number format/fixed},
log ticks with fixed point,
width=0.951\figurewidth,
height=\figureheight,
at={(0\figurewidth,0\figureheight)},
scale only axis,
xmin=50000,
xmax=125000,
xminorticks=true,
xtick={ 50000, 75000, 100000, 125000},
xlabel style={font=\color{white!15!black}},
xlabel={Matrix Size n},
ymode=log,
ymin=100,
ymax=20000,
yminorticks=true,
ylabel style={font=\color{white!15!black}},
ylabel={Runtime in s},
axis background/.style={fill=white},
axis x line*=bottom,
axis y line*=left,
xmajorgrids,
xminorgrids,
ymajorgrids,
yminorgrids,
legend columns=2,
legend style={at={(0.5,-.2)}, anchor=north, legend cell align=left, align=left, draw=white!15!black, font=\footnotesize, draw=none},
every axis plot/.append style={thick}
]
\addplot [color=black, dashed, mark=otimes, mark options={solid, black}]
  table[row sep=crcr]{%
50000	555.150448\\
75000	1892.597726\\
75000	1883.990606\\
100000	3887.45453\\
125000	7630.343413\\
};
\addlegendentry{Complex ELPA1, 100\%}

\addplot [color=blue, dashed, mark=triangle, mark options={solid, blue}]
  table[row sep=crcr]{%
50000	473.833636\\
75000	1624.282363\\
75000	1618.400776\\
100000	3272.951259\\
125000	6409.983405\\
};
\addlegendentry{Complex ELPA1, 50\%}

\addplot [color=black, mark=otimes, mark options={solid, black}]
  table[row sep=crcr]{%
50000	582.523669\\
75000	1873.634627\\
100000	4314.206719\\
125000	8319.805634\\
};
\addlegendentry{Complex ELPA2, 100\%}

\addplot [color=blue, mark=triangle, mark options={solid, blue}]
  table[row sep=crcr]{%
50000	381.112962\\
75000	1197.961772\\
100000	2718.745813\\
125000	5220.239052\\
};
\addlegendentry{Complex ELPA2, 50\%}

\addplot [color=black, dotted, mark=x, mark options={solid, black}]
  table[row sep=crcr]{%
50000	680.660579\\
75000	2170.239685\\
100000	5041.568116\\
125000	9707.702009\\
};
\addlegendentry{Complex MKL 100\%}

\addplot [color=blue, dotted, mark=+, mark options={solid, blue}]
  table[row sep=crcr]{%
50000	550.734223\\
75000	1752.525692\\
100000	4002.68405\\
125000	7803.777104\\
};
\addlegendentry{Complex MKL 50\%}

\addplot [color=red, dashed, mark=o, mark options={solid, red}]
  table[row sep=crcr]{%
50000	276.581136\\
75000	905.095688\\
75000	905.460716\\
100000	1984.195268\\
100000	1976.866348\\
125000	3749.649961\\
};
\addlegendentry{Skew-Symmetric\\ ELPA1, 100\%}

\addplot [color=black!50!green, dashed, mark=square, mark options={solid, black!50!green}]
  table[row sep=crcr]{%
50000	230.961236\\
75000	752.11202\\
75000	754.362344\\
100000	1597.896336\\
100000	1592.612259\\
125000	3015.884706\\
};
\addlegendentry{Skew-Symmetric\\ ELPA1, 50\%}

\addplot [color=red, mark=o, mark options={solid, red}]
  table[row sep=crcr]{%
50000	264.599098\\
75000	827.439294\\
100000	1884.332301\\
125000	3607.159364\\
150000	6037.844053\\
200000	13988.766414\\
};
\addlegendentry{Skew-Symmetric\\ ELPA2, 100\%}


\addplot [color=black!50!green, mark=square, mark options={solid, black!50!green}]
  table[row sep=crcr]{%
50000	164.182404\\
75000	502.653266\\
75000	498.541584\\
100000	1122.360256\\
100000	1120.392337\\
125000	2125.337833\\
};
\addlegendentry{Skew-Symmetric\\ ELPA2, 50\%}

\end{axis}
\end{tikzpicture}%

%% file: 1NodeGPU.tex
%
%
\definecolor{mycolor1}{rgb}{0.00000,0.44700,0.74100}%
\definecolor{mycolor2}{rgb}{0.85000,0.32500,0.09800}%
\begin{tikzpicture}

\begin{axis}[%
/pgf/number format/.cd,
        1000 sep={\,},
scaled x ticks = false,
x tick label style = {/pgf/number format/fixed},
log ticks with fixed point,
width=0.951\figurewidth,
height=\figureheight,
at={(0\figurewidth,0\figureheight)},
scale only axis,
xmin=0,
xmax=32768,
xminorticks=true,
xtick={ 1024,8192,16384,24576,32768 },
xlabel style={font=\color{white!15!black}},
xlabel={Matrix Size n},
ymin=0,
ymax=600,
yminorticks=true,
ylabel style={font=\color{white!15!black}},
ylabel={Runtime in s},
axis background/.style={fill=white},
axis x line*=bottom,
axis y line*=left,
xmajorgrids,
ymajorgrids,
legend columns=1,
legend style={at={(0.5,-.2)}, anchor=north, legend cell align=left, align=left, draw=white!15!black, font=\footnotesize, draw=none},
every axis plot/.append style={thick}
]
\addplot [color=mycolor1, mark=x, mark options={solid, mycolor1}]
  table[row sep=crcr]{%
1024	0.096474\\
2048	0.249621\\
4096	1.40249\\
8192	12.607739\\
16384	91.031407\\
24576	283.636567\\
32768	569.981855\\
};
\addlegendentry{ 2x Intel Xeon Silver 4110, ELPA2}
\addplot [color=mycolor2, mark=o, mark options={solid, mycolor2}]
  table[row sep=crcr]{%
1024	4.306552\\
2048	7.022835\\
4096	12.755546\\
8192	27.295248\\
16384	78.019717\\
24576	173.052262\\
32768	327.612085\\
};
\addlegendentry{ 2x Intel Xeon Silver 4110 + 1x Nvidia P100, ELPA1}
\end{axis}
\end{tikzpicture}%

%% file: AccBSE.tex
%
%
\definecolor{mycolor1}{rgb}{0.00000,0.44700,0.74100}%
\definecolor{mycolor2}{rgb}{0.85000,0.32500,0.09800}%
\definecolor{mycolor3}{rgb}{0.92900,0.69400,0.12500}%
\definecolor{mycolor4}{rgb}{0.49400,0.18400,0.55600}%
\begin{tikzpicture}

\begin{axis}[%
/pgf/number format/.cd,
        1000 sep={\,},
log basis x=16,
log ticks with fixed point,
xtick={ 16,  32,  64, 128, 256, 512},
xlabel={Number of cores},
width=0.951\figurewidth,
height=\figureheight,
at={(0\figurewidth,0\figureheight)},
scale only axis,
xmode=log,
xmin=16,
xmax=512,
ylabel={Runtime in $s$},
xminorticks=true,
ymode=log,
ymin=100,
ymax=20000,
yminorticks=true,
axis background/.style={fill=white},
axis x line*=bottom,
axis y line*=left,
xmajorgrids,
xminorgrids,
ymajorgrids,
yminorgrids,
legend columns=2,
legend style={at={(0.5,-.2)}, anchor=north, legend cell align=left, align=left, draw=none, font=\footnotesize},
every axis plot/.append style={thick}
]
\addplot [color=mycolor3, dashed, mark=square*, mark options={solid, mycolor3, fill}]
  table[row sep=crcr, stack plots=false]{%
16	10351.7796609402\\
32	6758.68391108513\\
64	5060.25346708298\\
128	3136.28374218941\\
256	1859.97967100143\\
512	1259.60617399216\\
};
\addlegendentry{BSEPACK NB = 64}

\addplot [color=mycolor4, dashed, mark=square*, mark options={solid, mycolor4, fill}]
  table[row sep=crcr, stack plots=false]{%
16	5827.98025894165\\
32	3270.47756314278\\
64	1748.30821108818\\
128	1014.14898109436\\
256	579.556116104126\\
512	377.349285125732\\
};
\addlegendentry{BSEPACK NB = 256}

\addplot[draw=black, mark=x, mark size=3pt] table[row sep=crcr]{%
16	2700.33380913734\\
32	1408.83185696602\\
64	735.923701047897\\
128	402.60865688324\\
256	225.283646821976\\
512	136.56920003891\\
};
\addlegendentry{BSEPACK + ELPA2}

\end{axis}
\end{tikzpicture}%